\documentstyle[amssymb,amsmath,graphicx,12pt]{article}
\newtheorem{th}{Theorem}[section]
\newtheorem{lem}[th]{Lemma}

\newtheorem{defn}[th]{Definition}
\newenvironment{defn-new}{\begin{defn} \em}{\end{defn}}
\newtheorem{rem}[th]{Remark}
\newenvironment{rem-new}{\begin{rem} \em}{\end{rem}}
\newtheorem{ex}[th]{Example}
\newenvironment{ex-new}{\begin{ex} \em}{\end{ex}}
\newtheorem{prob}[th]{Problem}
\newenvironment{prob-new}{\begin{prob} \em}{\end{prob}}

\newenvironment{notation-new}{\begin{rem} \em}{\end{rem}}

\newenvironment{agr-new}{\begin{rem} \em}{\end{rem}}

\makeatletter \@addtoreset{equation}{section} \makeatother
\begin{document}


\begin{center}
{\LARGE {\bf Inequalities for algebraic Casorati curvatures and their
applications}}
\end{center}

\medskip

\begin{center}
{\em Dedicated to Felice Casorati}
\end{center}

\bigskip

\begin{center}
Mukut Mani Tripathi\footnote{Submitted to Note di Mat}
\end{center}

\bigskip


\begin{quote}
\noindent {\bf Abstract.} The notion of different kind of algebraic
Casorati curvatures are introduced. Some results expressing basic
Casorati inequalities for algebraic Casorati curvatures are
presented. Equality cases are also discussed. As a simple
application, basic Casorati inequalities for different $\delta
$-Casorati curvatures for Riemannian submanifolds are presented.
Further applying these results, Casorati inequalities for Riemannian
submanifolds of real space forms are obtained. Finally, some problems
are presented for further studies.
\end{quote}


\begin{quote}
\noindent {\em AMS $2000$ Mathematics Subject Classification\/}: 53B20.
\end{quote}


\begin{quote}
\end{quote}

%
%
%

\section{Introduction\label{sect-intro}}

Felice Casorati was one of the great Italian mathematicians. best
known for the Casorati-Weierstrass theorem in complex analysis. He
was born in Pavia on December 17, 1835 and his soul departed on
September 11, 1890 in Casteggio. Before his departure, in 1889,
Casorati \cite{Casorati-1889} defined a curvature for a regular
surface in Euclidean {$3$-space} which turns out to be the
normalized sum of the squared principal curvatures. In
\cite{Casorati-1889-note}, the author says that he could not check
the paper \cite{Casorati-1889} before printing, and advices readers
to rather use a subsequent paper \cite{Casorati-1890}. This
curvature is now well known as the Casorati curvature. Several
geometers believe that Casorati preferred this curvature over the
traditional Gaussian curvature because the Casorati curvature
vanishes for a surface in Euclidean $3$-space if and only if both
Euler normal curvatures (or principal curvatures) of the surface
vanish simultaneously and thus corresponds better with the common
intuition of curvature. For a hypersurface of a Riemannian manifold
the Casorati curvature is defined to be the normalized sum of the
squared principal normal curvatures of the hypersurface, and in
general, the Casorati curvature of a submanifold of a Riemannian
manifold is defined to
be the normalized squared length of the second fundamental form \cite%
{Decu-HV-08-JIPAM}. Geometrical meaning and the importance of the Casorati
curvature, discussed by several geometers, can be visualized in several
research/survey papers including \cite{Decu-08-Brasov}, \cite%
{Decu-PPV-13-Krag}, \cite{Decu-PSV-14}, \cite{Ghiso-08-Brasov}, \cite%
{Haesen-KV-09-Note-Mat}, \cite{Koenderink-13-handbook}, \cite%
{Kowalczyk-08-Brasov}, \cite{Ons-Vert-11-Krag}, \cite{Verst-04-Soochow} and
\cite{Verst-13-Krag}.

The paper is organized as follows. In section~\ref%
{sect-curvature-like-tensor}, some preliminaries regarding curvature like
tensors are presented. In section~\ref{sect-Casorati-alg-Casorati-curvatures}%
, given an $n$-dimensional Riemannian manifold $\left( M,g\right) $, a
Riemannian vector bundle $\left( B,g_{B}\right) $ over $M$, a $B$-valued
symmetric $\left( 1,2\right) $-tensor field $\zeta $ and a (curvature-like)
tensor field $T$ satisfying the algebraic Gauss equation, we introduce the
notion of different kind of algebraic Casorati curvatures $\widehat{\delta }%
_{{\cal C}^{T,\zeta }}(n-1)$, $\delta _{{\cal C}^{T,\zeta }}(n-1)$, $\delta
_{{\cal C}^{T,\zeta }}(r;n-1)$, $\widehat{\delta }_{{\cal C}^{T,\zeta
}}(r;n-1)$, which in special cases of Riemannian submanifolds reduce to
already known $\delta $-Casorati curvatures. In section~\ref%
{sect-Casorati-ineq-alg}, first we prove an useful Lemma regarding a
constrained extremum problem. Then we present results expressing basic
Casorati inequalities for algebraic Casorati curvatures. Equality cases are
also discussed. After this, application parts begin. In section~\ref%
{sect-Casorati-ineq-Riem-submfd}, we obtain basic Casorati inequalities for
Casorati curvatures $\delta (r;n-1)$, $\widehat{\delta }(r;n-1)$, $\delta
(n-1)$, $\widehat{\delta }(n-1)$ for Riemannian submanifolds. In section~\ref%
{sect-Casorati-ineq-Riem-submfd-RSF}, we further apply these results to
obtain Casorati inequalities for Riemannian submanifolds of real space forms
with very short proofs. Finally, in section~\ref%
{sect-Casorati-further-studies}, we present some problems for further
studies.

\section{Curvature like tensor\label{sect-curvature-like-tensor}}

In $1967$, R.S. Kulkarni introduced the notion of a {\em curvature structure}
(cf. \cite[{\S }$8$ of Chapter 1]{Kulk-67-thesis}, \cite{Kulk-69-BAMS}),
which is now widely known as a {\em curvature-like tensor} ({\em field}).
Let $(M,g)$ be an {$n$-dimensional} Riemannian manifold. Let $T$ be a
curvature-like tensor so that it satisfies the following properties
\begin{equation}
T(X,Y,Z,W)=-\,T(Y,X,Z,W),  \label{eq-CLT-1}
\end{equation}%
\begin{equation}
T(X,Y,Z,W)=T(Z,W,X,Y),  \label{eq-CLT-3}
\end{equation}%
\begin{equation}
T(X,Y,Z,W)+T(Y,Z,X,W)+T(Z,X,Y,W)=0  \label{eq-CLT-5}
\end{equation}%
for all vector fields $X$, $Y$, $Z$ and $W$ on $M$. For a curvature-like
tensor field $T$, the $T${\em -sectional curvature} associated with a $2$%
-plane section $\Pi _{2}$ spanned by orthonormal vectors $X$ and $Y$ at $%
p\in M$, is given by \cite{Bolt-DFV-09-MIA}
\[
K_{T}(\Pi _{2})=K_{T}(X\wedge Y)=T(X,Y,Y,X).
\]%
Let $\left\{ e_{1},e_{2},\ldots ,e_{n}\right\} $ be any orthonormal basis of
$T_{p}M$. The $T${\em -Ricci tensor} $\,S_{T}$ is defined by
\[
S_{T}(X,Y)=\sum_{j=1}^{n}T\left( e_{j},X,Y,e_{j}\right) ,\qquad X,Y\in
T_{p}M.
\]%
The $T${\em -Ricci curvature} is given by
\[
{\rm Ric}_{T}(X)=S_{T}(X,X),\qquad X\in T_{p}M.
\]%
The $T${\em -scalar curvature} is given by \cite{Bolt-DFV-09-MIA}
\begin{equation}
\tau _{T}(p)=\sum_{1\leq i<j\leq n}T\left( e_{i},e_{j},e_{j},e_{i}\right) ,
\label{eq-tau-T-p}
\end{equation}%
Now, let $\Pi _{k}$ be a $k$-plane section of $T_{p}M$ and $X$ a unit vector
in $\Pi _{k}$. If $k=n$ then $\Pi _{n}=T_{p}M$; and if $k=2$ then $\Pi _{2}$
is a plane section of $T_{p}M$. We choose an orthonormal basis $%
\{e_{1},\ldots ,e_{k}\}$ of $\Pi _{k}$. Then we define the $T$-$k${\em %
-Ricci curvature} of $\Pi _{k}$ at $e_{i}$, $i\in \{1,\ldots ,k\}$, denoted $%
({\rm Ric}_{T})_{\Pi _{k}}(e_{i})$, by
\begin{equation}
({\rm Ric}_{T})_{\Pi _{k}}(e_{i})=\sum_{j=1,j\neq i}^{k}K_{T}(e_{i}\wedge
e_{j}).  \label{eq-Ric-Pi-k-ei}
\end{equation}%
We note that a {$T$-$n$-{\em Ricci curvature}} $({\rm Ric}%
_{T})_{T_{p}M}(e_{i})$ is the usual $T${\em -Ricci curvature} of $e_{i}$,
denoted ${\rm Ric}_{T}(e_{i})$. The $T${\em -}$k${\em -scalar curvature} $%
\tau _{T}(\Pi _{k})$ of the $k$-plane section $\Pi _{k}$ is given by
\begin{equation}
\tau _{T}(\Pi _{k})=\sum_{1\leq i<j\leq k}K_{T}(e_{i}\wedge e_{j}).
\label{eq-tau-T-Pi-k-1}
\end{equation}%
We note that
\begin{equation}
\tau _{T}(\Pi _{k})=\frac{1}{2}\,\sum_{i=1}^{k}\sum_{j=1,j\neq
i}^{k}K_{T}(e_{i}\wedge e_{j})=\frac{1}{2}\,\sum_{i=1}^{k}({\rm Ric}%
_{T})_{\Pi _{k}}(e_{i}).  \label{eq-tau-T-Pi-k-2}
\end{equation}%
The $T$-scalar curvature of $M$ at $p$ is identical with the $T$-$n$-scalar
curvature of the tangent space $T_{p}M$ of $M$ at $p$, that is, $\tau
_{T}(p)=\tau _{T}(T_{p}M)$. If $\Pi _{2}$ is a $2$-plane section, $\tau
_{T}(\Pi _{2})$ is nothing but the $T$-sectional curvature $K_{T}(\Pi _{2})$
of $\Pi _{2}$. The $T${\em -}$k${\em -normalized scalar curvature} of a $k$%
-plane section $\Pi _{k}$ at $p$ is defined as
\[
(\tau _{T})_{{\rm Nor}}(\Pi _{k})=\frac{2}{k(k-1)}\,\tau _{T}(\Pi _{k}).
\]%
The $T${\em -normalized scalar curvature} at $p$ is defined as
\[
(\tau _{T})_{{\rm Nor}}(p)=(\tau _{T})_{{\rm Nor}}(T_{p}M)=\frac{2}{n(n-1)}%
\,\tau _{T}(p).
\]

If $T$ is replaced by the Riemann curvature tensor $R$, then $T$-sectional
curvature $K_{T}$, $T$-Ricci tensor $S_{T}$, $T$-Ricci curvature ${\rm Ric}%
_{T}$, $T$-scalar curvature $\tau _{T}$, $T$-normalized scalar curvature $%
(\tau _{T})_{{\rm Nor}}$, $T$-$k$-Ricci curvature $({\rm Ric}_{T})_{\Pi
_{k}} $, $T$-$k$-scalar curvature $\tau _{T}(\Pi _{k})$, $T$-$k$-normalized
scalar curvature $(\tau _{T})_{{\rm Nor}}(\Pi _{k})$ and $T$-normalized
scalar curvature $(\tau _{T})_{{\rm Nor}}$ become the sectional curvature $K$%
, the Ricci tensor $S$, the Ricci curvature ${\rm Ric}$, the scalar
curvature $\tau $, the normalized scalar curvature $\tau _{{\rm Nor}}$, $k$%
-Ricci curvature ${\rm Ric}_{\Pi _{k}}$, $k$-scalar curvature $\tau (\Pi
_{k})$, $k$-normalized scalar curvature $\tau _{{\rm Nor}}(\Pi _{k})$ and
normalized scalar curvature $\tau _{{\rm Nor}}$, respectively.

\section{Algebraic Casorati curvatures\label%
{sect-Casorati-alg-Casorati-curvatures}}

Let $\left( M,g\right) $ be an {$n$-dimensional} submanifold of an {$m$%
-dimensional} Riemannian manifold $(\widetilde{M},\widetilde{g})$. The
equation of Gauss is given by
\begin{eqnarray}
R(X,Y,Z,W) &=&\widetilde{R}(X,Y,Z,W)+\ \widetilde{g}\left( \sigma
(Y,Z),\sigma (X,W)\right)  \nonumber \\
&&-\ \widetilde{g}\left( \sigma (X,Z),\sigma (Y,W)\right)  \label{eq-Gauss}
\end{eqnarray}%
for all $X,Y,Z,W\in TM$, where $\widetilde{R}$ and $R$ are the curvature
tensors of $\widetilde{M}$ and $M$, respectively and $\sigma $ is the second
fundamental form of the immersion of $M$ in $\widetilde{M}$. The Ricci-K\"{u}%
hn equation is given by
\begin{equation}
R^{\perp }(X,Y,N,V)=\widetilde{R}(X,Y,N,V)+g\left( \left[ A_{N},A_{V}\right]
X,Y\right)  \label{eq-Ricci-Kuhn}
\end{equation}%
for all $X,Y\in TM$ and for all $N,V\in T^{\perp }M$, where
\[
R^{\perp }(X,Y)N=\nabla _{X}^{\perp }\nabla _{Y}^{\perp }N-\nabla
_{Y}^{\perp }\nabla _{X}^{\perp }N-\nabla _{\left[ X,Y\right] }^{\perp }N,
\]%
\[
\left[ A_{N},A_{V}\right] =A_{N}A_{V}-A_{V}A_{N},
\]%
with $\nabla ^{\perp }$ being the induced normal connection in the normal
bundle $T^{\perp }M$ and $A_{N}$ being the shape operator in the direction $%
N $.

Let $M$ be an {$n$-dimensional} Riemannian submanifold of an {$m$-dimensional%
} Riemannian manifold $\widetilde{M}$. A point $p\in M$ is said to be an
{\em invariantly quasi-umbilical point} if there exist $m-n$ mutually
orthogonal unit normal vectors $N_{n+1},\ldots ,N_{m}$ such that the shape
operators with respect to all directions $N_{\alpha }$ have an eigenvalue of
multiplicity $n-1$ and that for each $N_{\alpha }$ the distinguished
eigendirection is the same. The submanifold is said to be an {\em %
invariantly quasi-umbilical submanifold} if each of its points is an
invariantly quasi-umbilical point. For details, we refer to \cite%
{Bl-Led-77-Stevin}.

Let $\left( M,g\right) $ be an {$n$-dimensional} Riemannian submanifold of
an {$m$-dimensional} Riemannian manifold $(\widetilde{M},\widetilde{g})$.
Let $\left\{ e_{1},\ldots ,e_{n}\right\} $ be an orthonormal basis of the
tangent space $T_{p}M$ and $e_{\alpha }$ belongs to an orthonormal basis $%
\left\{ e_{n+1},\ldots ,e_{m}\right\} $ of the normal space $T_{p}^{\perp }M$%
. We let
\[
\sigma _{ij}^{\alpha }=\widetilde{g}\left( \sigma \left( e_{i},e_{j}\right)
,e_{\alpha }\right) ,\quad i,j\in \{1,\ldots ,n\},\quad \alpha \in
\{n+1,\ldots ,m\}.
\]%
Then, the squared mean curvature of the submanifold $M$ in $\widetilde{M}$
is defined by
\[
\left\Vert H\right\Vert ^{2}=\frac{1}{n^{2}}\sum_{\alpha =n+1}^{m}\left(
\sum_{i=1}^{n}\sigma _{ii}^{\alpha }\right) ^{2},
\]%
and the squared norm of second fundamental form $\sigma $ is
\[
\left\Vert \sigma \right\Vert ^{2}=\sum_{i,j=1}^{n}\widetilde{g}\left(
\sigma \left( e_{i},e_{j}\right) ,\sigma \left( e_{i},e_{j}\right) \right) .
\]%
Let $K_{ij}$ and $\widetilde{K}_{ij}$ denote the sectional curvature of the
plane section spanned by $e_{i}$ and $e_{j}$ at $p$ in the submanifold $M$
and in the ambient manifold $\widetilde{M}$, respectively. In view of (\ref%
{eq-Gauss}), we have
\begin{equation}
K_{ij}=\widetilde{K}_{ij}+\sum_{\alpha =n+1}^{m}\left( \sigma _{ii}^{\alpha
}\sigma _{jj}^{\alpha }-(\sigma _{ij}^{\alpha })^{2}\right) .  \label{eq-Kij}
\end{equation}%
From (\ref{eq-Kij}) it follows that
\begin{equation}
2\tau (p)=2\widetilde{\tau }\left( T_{p}M\right) +n^{2}\left\Vert
H\right\Vert ^{2}-\left\Vert \sigma \right\Vert ^{2},  \label{eq-tau-H-sigma}
\end{equation}%
where
\[
\widetilde{\tau }\left( T_{p}M\right) =\sum_{1\leq i<j\leq n}\widetilde{K}%
_{ij}
\]%
denotes the scalar curvature of the {$n$-plane} section $T_{p}M$ in the
ambient manifold $\widetilde{M}$. From (\ref{eq-tau-H-sigma}), it
immediately follows that
\begin{equation}
\tau _{{\rm Nor}}(p)=\widetilde{\tau }_{{\rm Nor}}\left( T_{p}M\right) +%
\frac{n}{n-1}\left\Vert H\right\Vert ^{2}-\frac{1}{n(n-1)}\left\Vert \sigma
\right\Vert ^{2},  \label{eq-tau-Nor-H-sigma}
\end{equation}%
where
\begin{equation}
\tau _{{\rm Nor}}(p)=\frac{2\tau (p)}{n(n-1)},\qquad \widetilde{\tau }_{{\rm %
Nor}}\left( T_{p}M\right) =\frac{2\widetilde{\tau }\left( T_{p}M\right) }{%
n(n-1)}.  \label{eq-tau-Nor-H-sigma-a}
\end{equation}

The {\em Casorati curvature} ${\cal C}$ \cite{Decu-HV-08-JIPAM} of the
Riemannian submanifold $M$ is defined to be the normalized squared length of
the second fundamental form $\sigma $, that is,
\begin{equation}
{\cal C}=\frac{1}{n}\left\Vert \sigma \right\Vert ^{2}=\frac{1}{n}%
\sum_{\alpha =n+1}^{m}\sum_{i,j=1}^{n}\left( \sigma _{ij}^{\alpha }\right)
^{2}.  \label{eq-C}
\end{equation}%
For a $k$-dimensional subspace $\Pi _{k}$ of $T_{p}M$, $k\geq 2$ spanned by $%
\{e_{1},\ldots ,e_{k}\}$, the Casorati curvature ${\cal C}\left( \Pi
_{k}\right) $ of the subspace $\Pi _{k}$ is defined to be \cite%
{Decu-HV-07-Brasov}
\[
{\cal C}\left( \Pi _{k}\right) =\frac{1}{k}\sum_{\alpha
=n+1}^{m}\sum_{i,j=1}^{k}\left( \sigma _{ij}^{\alpha }\right) ^{2}.
\]%
The {\em normalized} $\delta ${\em -Casorati curvatures} $\widehat{\delta }_{%
{\cal C}}(n-1)$, $\delta _{{\cal C}}^{\prime }(n-1)$ of a Riemannian
submanifold $M$ are given by \cite{Decu-HV-07-Brasov}
\begin{equation}
\lbrack \widehat{\delta }_{{\cal C}}(n-1)]_{p}=2\,{\cal C}_{p}-\frac{2n-1}{2n%
}\sup \left\{ {\cal C}(\Pi _{n-1}):\Pi _{n-1}\ {\rm is}\ {\rm a}\ {\rm %
hyperplane}\ {\rm of}\ T_{p}M\right\} ,  \label{eq-delta-hat-C-(n-1)}
\end{equation}

\begin{equation}
\lbrack \delta _{{\cal C}}^{\prime }(n-1)]_{p}=\frac{1}{2}\,{\cal C}_{p}+%
\frac{n+1}{2n(n-1)}\inf \left\{ {\cal C}(\Pi _{n-1}):\Pi _{n-1}\ {\rm is}\
{\rm a}\ {\rm hyperplane}\ {\rm of}\ T_{p}M\right\} .
\label{eq-delta-prime-C-(n-1)}
\end{equation}%
In \cite{Decu-HV-07-Brasov}, the authors denoted $\delta _{{\cal C}}^{\prime
}(n-1)$ by $\delta _{{\cal C}}(n-1)$. The (modified) {\em normalized} $%
\delta ${\em -Casorati curvatures} $\delta _{{\cal C}}(n-1)$ of the
Riemannian submanifold $M$ is given by (\cite{LeeCW-YL-14-JIA}, \cite%
{ZhangP-Z-15-arXiv})
\begin{equation}
\lbrack \delta _{{\cal C}}(n-1)]_{p}=\frac{1}{2}\,{\cal C}_{p}+\frac{n+1}{2n}%
\inf \left\{ {\cal C}(\Pi _{n-1}):\Pi _{n-1}\ {\rm is}\ {\rm a}\ {\rm %
hyperplane}\ {\rm of}\ T_{p}M\right\} .  \label{eq-delta-C-(n-1)}
\end{equation}%
It should be noted that the normalized {$\delta $-Casorati} curvatures $%
\widehat{\delta }_{{\cal C}}(n-1)$, $\delta _{{\cal C}}^{\prime }(n-1)$ and $%
\delta _{{\cal C}}(n-1)$ vanish trivially for $n=2$ \cite{ZhangP-Z-15-arXiv}%
. In \cite{LeeCW-YL-14-JIA}, the authors pointed out that the coefficient $%
\frac{n+1}{2n(n-1)}$ in (\ref{eq-delta-prime-C-(n-1)}) was inappropriate and
therefore they modified the coefficient from $\frac{n+1}{2n(n-1)}$ to $\frac{%
n+1}{2n}$ in the definition of $\delta _{{\cal C}}^{\prime }(n-1)$ to obtain
the definition of $\delta _{{\cal C}}(n-1)$ (see also \cite%
{LeeJW-V-15-Taiwan}). For a positive real number $r\neq n(n-1)$, letting
\[
{\rm a}(r)=\frac{1}{nr}(n-1)\left( n+r\right) \left( n^{2}-n-r\right) ,
\]%
the {\em normalized} $\delta $-{\em Casorati curvatures} $\delta _{{\cal C}%
}(r;n-1)$ and $\widehat{\delta }_{{\cal C}}(r;n-1)$ of a Riemannian
submanifold $M$ are given by \cite{Decu-HV-08-JIPAM}
\begin{equation}
\lbrack \delta _{{\cal C}}(r;n-1)]_{p}=r\,{\cal C}_{p}+{\rm a}(r)\inf
\left\{ {\cal C}(\Pi _{n-1}):\Pi _{n-1}\ {\rm is}\ {\rm a}\ {\rm hyperplane}%
\ {\rm of}\ T_{p}M\right\} ,  \label{eq-delta-C-(r;n-1)}
\end{equation}%
if $0<r<n(n-1)$, and
\begin{equation}
\lbrack \widehat{\delta }_{{\cal C}}(r;n-1)]_{p}=r\,{\cal C}_{p}+{\rm a}%
(r)\sup \left\{ {\cal C}(\Pi _{n-1}):\Pi _{n-1}\ {\rm is}\ {\rm a}\ {\rm %
hyperplane}\ {\rm of}\ T_{p}M\right\} ,  \label{eq-delta-hat-C-(r;n-1)}
\end{equation}%
if $n(n-1)<r$, respectively.

In \cite{LeeCW-YL-14-JIA} the normalized {$\delta $-Casorati} curvatures $%
\widehat{\delta }_{{\cal C}}(r;n-1)$ and $\delta _{{\cal C}}(r;n-1)$ are
called as the generalized normalized {$\delta $-Casorati} curvatures $%
\widehat{\delta }_{{\cal C}}(r;n-1)$ and $\delta _{{\cal C}}(r;n-1)$,
respectively. We see that \cite{LeeJW-V-15-Taiwan}
\begin{equation}
\left[ \delta _{{\cal C}}(n-1)\right] _{p}=\frac{1}{n(n-1)}\left[ \delta _{%
{\cal C}}\left( \frac{n(n-1)}{2};n-1\right) \right] _{p},
\label{eq-delta-C-(r;n-1)-C-(n-1)}
\end{equation}%
\begin{equation}
\left[ \widehat{\delta }_{{\cal C}}(n-1)\right] _{p}=\frac{1}{n(n-1)}\left[
\widehat{\delta }_{{\cal C}}\left( 2n(n-1);n-1\right) \right] _{p}
\label{eq-delta-hat-C-(r;n-1)-C-(n-1)}
\end{equation}%
for all $p\in M$.

Let $\left( M,g\right) $ be an {$n$-dimensional} Riemannian manifold and $%
\left( B,g_{B}\right) $ a Riemannian vector bundle over $M$. If $\zeta $ is
a $B$-valued symmetric $\left( 1,2\right) $-tensor field and $T$ a $\left(
0,4\right) $-tensor field on $M$ such that
\begin{equation}
T(X,Y,Z,W)=g_{B}(\zeta (X,W),\zeta (Y,Z))-g_{B}(\zeta (X,Z),\zeta (Y,W))
\label{eq-Gauss-alg}
\end{equation}%
for all vector fields $X$,$Y$,$Z$,$W$ on $M$, then the equation (\ref%
{eq-Gauss-alg}) is said to be an {\em algebraic Gauss equation} \cite%
{Chen-DV-05-DGA}. Every $(0,4)$-tensor field $T$ on $M$, which satisfies (%
\ref{eq-Gauss-alg}), becomes a curvature-like tensor.

A typical example of an algebraic Gauss equation is given for a submanifold $%
M$ of an Euclidean space, if $B$ is the normal bundle, $\zeta $ the second
fundamental form and $T$ the curvature tensor. Some nice situations, in
which such $T$ and $\zeta $ satisfying an algebraic Gauss equation exist,
are Lagrangian and Kaehlerian slant submanifolds of complex space forms and $%
C$-totally real submanifolds of Sasakian space forms.

Now, let $\{e_{1},\ldots ,e_{n}\}$ be an orthonormal basis of the tangent
space $T_{p}M$ and $e_{\alpha }$ belong to an orthonormal basis $%
\{e_{n+1},\ldots ,e_{m}\}$ of the Riemannian vector bundle $(B,g_{B})$ over $%
M$ at $p$. We put
\[
\zeta _{ij}^{\alpha }=g_{B}\left( \zeta \left( e_{i},e_{j}\right) ,e_{\alpha
}\right) ,\quad \left\Vert \zeta \right\Vert
^{2}=\sum_{i,j=1}^{n}g_{B}\left( \zeta \left( e_{i},e_{j}\right) ,\zeta
\left( e_{i},e_{j}\right) \right) ,
\]%
\[
{\rm trace}\,\zeta =\sum_{i=1}^{n}\zeta \left( e_{i},e_{i}\right) ,\quad
\left\Vert {\rm trace}\,\zeta \right\Vert ^{2}=g_{B}({\rm trace}\,\zeta ,%
{\rm trace}\,\zeta ).
\]

Motivated by the definitions given in \cite{Decu-HV-07-Brasov}, \cite%
{Decu-HV-08-JIPAM} and \cite{LeeCW-YL-14-JIA} we give the following
definitions.

\begin{defn-new}
\label{defn-alg-Casorati-curvature-1} Let $\left( M,g\right) $ be an {$n$%
-dimensional} Riemannian manifold, $\left( B,g_{B}\right) $ a Riemannian
vector bundle over $M$, $\zeta $ a {$B$-valued} symmetric {$\left(
1,2\right) $-tensor} field on $M$, and $T$ a curvature-like tensor field
satisfying the algebraic Gauss equation $(\ref{eq-Gauss-alg})$. Then the
{\em algebraic Casorati curvature} ${\cal C}^{T,\zeta }$ with respect to $T$
and the Riemannian vector bundle $\left( B,g_{B}\right) $ over $M$ is
defined to be
\begin{equation}
{\cal C}^{T,\zeta }=\frac{1}{n}\left\Vert \zeta \right\Vert ^{2}=\frac{1}{n}%
\sum_{\alpha =n+1}^{m}\sum_{i,j=1}^{n}\left( \zeta _{ij}^{\alpha }\right)
^{2}.  \label{eq-alg-Casorati-curv}
\end{equation}%
For a $k$-dimensional subspace $\Pi _{k}$ of $T_{p}M$, $k\geq 2$, spanned by
$\{e_{1},\ldots ,e_{k}\}$, the {\em algebraic Casorati curvature} ${\cal C}%
^{T,\zeta }(\Pi _{k})$ of the subspace $\Pi _{k}$ is defined to be
\begin{equation}
{\cal C}^{T,\zeta }(\Pi _{k})=\frac{1}{k}\sum_{\alpha
=n+1}^{m}\sum_{i,j=1}^{k}\left( \zeta _{ij}^{\alpha }\right) ^{2}.
\label{eq-alg-Casorati-curv-Pi-k}
\end{equation}
\end{defn-new}

We note that
\[
{\cal C}_{p}^{T,\zeta }={\cal C}^{T,\zeta }(T_{p}M),\qquad p\in M.
\]

\begin{defn-new}
\label{defn-alg-Casorati-curvature-2} Let $\left( M,g\right) $ be an {$n$%
-dimensional} Riemannian manifold, $\left( B,g_{B}\right) $ a Riemannian
vector bundle over $M$, $\zeta $ a $B$-valued symmetric $\left( 1,2\right) $%
-tensor field on $M$, and $T$ a curvature-like tensor field satisfying the
algebraic Gauss equation {\rm (\ref{eq-Gauss-alg})}. Then we define the
following three {\em algebraic Casorati curvatures} $\delta _{{\cal C}%
^{T,\zeta }}(n-1)$ and $\widehat{\delta }_{{\cal C}^{T,\zeta }}(n-1)$ and by
\begin{eqnarray}
\lbrack \delta _{{\cal C}^{T,\zeta }}(n-1)]_{p} &=&\frac{1}{2}\,{\cal C}%
_{p}^{T,\zeta }+\frac{n+1}{2n}\inf \left\{ {\cal C}^{T,\zeta }(\Pi
_{n-1}):\Pi _{n-1}\ {\rm is}\ {\rm a}\ {\rm hyperplane}\ {\rm of}\
T_{p}M\right\} ,  \nonumber \\
&&  \label{eq-delta-C-(T,B)-(n-1)}
\end{eqnarray}%
\begin{eqnarray}
\lbrack \widehat{\delta }_{{\cal C}^{T,\zeta }}(n-1)]_{p} &=&2\,{\cal C}%
_{p}^{T,\zeta }-\frac{2n-1}{2n}\sup \{{\cal C}^{T,\zeta }(\Pi _{n-1}):\Pi
_{n-1}\ {\rm is}\ {\rm a}\ {\rm hyperplane}\ {\rm of}\ T_{p}M\}.  \nonumber
\\
&&  \label{eq-delta-hat-C-(T,B)-(n-1)}
\end{eqnarray}
\end{defn-new}

\begin{defn-new}
\label{defn-alg-Casorati-curvature-3} Let $\left( M,g\right) $ be an {$n$%
-dimensional} Riemannian manifold, $\left( B,g_{B}\right) $ a Riemannian
vector bundle over $M$, $\zeta $ a $B$-valued symmetric $\left( 1,2\right) $%
-tensor field on $M$, and $T$ a curvature-like tensor field satisfying the
algebraic Gauss equation {\rm (\ref{eq-Gauss-alg})}. For a positive real
number $r\neq n(n-1)$, let
\[
{\rm a}(r)=\frac{1}{nr}(n-1)\left( n+r\right) \left( n^{2}-n-r\right)
\]%
and define the {\em algebraic Casorati curvatures} $\delta _{{\cal C}%
^{T,\zeta }}(r;n-1)$ and $\widehat{\delta }_{{\cal C}^{T,\zeta }}(r;n-1)$ by
\begin{eqnarray}
\lbrack \delta _{{\cal C}^{T,\zeta }}(r;n-1)]_{p} &=&r\,{\cal C}%
_{p}^{T,\zeta }+{\rm a}(r)\inf \left\{ {\cal C}^{T,\zeta }(\Pi _{n-1}):\Pi
_{n-1}\ {\rm is}\ {\rm a}\ {\rm hyperplane}\ {\rm of}\ T_{p}M\right\}
\nonumber \\
&&  \label{eq-delta-C-(T,B)-(r;n-1)}
\end{eqnarray}%
\noindent if $0<r<n(n-1)$, and
\begin{eqnarray}
\lbrack \widehat{\delta }_{{\cal C}^{T,\zeta }}(r;n-1)]_{p} &=&r\,{\cal C}%
_{p}^{T,\zeta }+{\rm a}(r)\sup \left\{ {\cal C}^{T,\zeta }(\Pi _{n-1}):\Pi
_{n-1}\ {\rm is}\ {\rm a}\ {\rm hyperplane}\ {\rm of}\ T_{p}M\right\}
\nonumber \\
&&  \label{eq-delta-C-hat-(T,B)-(r;n-1)}
\end{eqnarray}%
\noindent if $n(n-1)<r$.
\end{defn-new}

\begin{rem-new}
\label{rem-alg-Casorati-curvature-1} Let $\left( M,g\right) $ be an {$n$%
-dimensional} Riemannian submanifold of an {$m$-dimensional} Riemannian
manifold $(\widetilde{M},\widetilde{g})$. Let the Riemannian vector bundle $%
\left( B,g_{B}\right) $ over $M$ be replaced by the normal bundle $T^{\perp
}M$, and the $B$-valued symmetric $\left( 1,2\right) $-tensor field $\zeta $
be replaced by the second fundamental form of immersion $\sigma $. Then the
algebraic Casorati curvature ${\cal C}^{T,\zeta }$ becomes the {\em Casorati
curvature} ${\cal C}$ of the Riemannian submanifold $M$ given by (\ref{eq-C}%
). 
The algebraic Casorati curvatures$\delta _{{\cal C}^{T,\zeta }}(n-1)$ and $%
\widehat{\delta }_{{\cal C}^{T,\zeta }}(n-1)$ become {\em normalized} $%
\delta ${\em -Casorati curvatures} $\delta _{{\cal C}}(n-1)$ and $\widehat{%
\delta }_{{\cal C}}(n-1)$ of the Riemannian submanifold $M$ given by $(\ref%
{eq-delta-C-(n-1)})$ and $(\ref{eq-delta-hat-C-(n-1)})$, respectively.
Finally, algebraic Casorati curvatures $\delta _{{\cal C}^{T,\zeta }}\left(
r;n-1\right) $ and $\widehat{\delta }_{{\cal C}^{T,\zeta }}\left(
r;n-1\right) $ become normalized {$\delta $-Casorati} curvatures $\delta _{%
{\cal C}}(r;n-1) $ and \allowbreak $\widehat{\delta }_{{\cal C}}(r;n-1) $ of
the Riemannian submanifold $M$ given by $(\ref{eq-delta-C-(r;n-1)})$ and $(%
\ref{eq-delta-hat-C-(r;n-1)})$, respectively.
\end{rem-new}

Now, we present the following useful Lemma.

\begin{lem}
\label{lem-C-(T,B)-trace-zeta-tau-T} Let $\left( M,g\right) $ be an {$n$%
-dimensional} Riemannian manifold, $\left( B,g_{B}\right) $ a Riemannian
vector bundle over $M$ and $\zeta $ a $B$-valued symmetric $\left(
1,2\right) $-tensor field. Let $T$ be a curvature-like tensor field
satisfying the algebraic Gauss equation {\rm (\ref{eq-Gauss-alg})}. Then
\begin{equation}
n{\cal C}^{T,\zeta }-\left\Vert {\rm trace}\,\zeta \right\Vert ^{2}=-\,2\tau
_{T}.  \label{eq-alg-nC-trace-zeta}
\end{equation}
\end{lem}

\noindent {\bf Proof.} Let $p\in M$, the set $\{e_{1},\ldots ,e_{n}\}$ be an
orthonormal basis of the tangent space $T_{p}M$ and $e_{\alpha }$ belong to
an orthonormal basis $\{e_{n+1},\ldots ,e_{m}\}$ of the Riemannian vector
bundle $(B,g_{B})$ over $M$ at $p$. From (\ref{eq-Gauss-alg}), we get
\begin{equation}
(K_{T})_{ij}=T\left( e_{i},e_{j},e_{j},e_{i}\right) =\sum_{\alpha
=n+1}^{m}\left( \zeta _{ii}^{\alpha }\zeta _{jj}^{\alpha }-(\zeta
_{ij}^{\alpha })^{2}\right) ,  \label{eq-alg-T-Kij}
\end{equation}%
which implies that
\begin{equation}
2\tau _{T}=\left\Vert {\rm trace}\,\zeta \right\Vert ^{2}-\Vert \zeta \Vert
^{2}=\left\Vert {\rm trace}\,\zeta \right\Vert ^{2}-n{\cal C}^{T,\zeta }.
\label{eq-alg-T-tau-trace-zeta-C}
\end{equation}%
This gives (\ref{eq-alg-nC-trace-zeta}). $\blacksquare $

\section{Basic Casorati inequalities\label{sect-Casorati-ineq-alg}}

We begin with the following two Lemmas:

\begin{lem}
\label{lem-optimization-convex-1} {\rm (\cite[Theorem~21.4, p.~425]%
{Chong-Zak-01-book})} Let $\Upsilon \subset {\Bbb R}^{n}$ be an open convex
set in ${\Bbb R}^{n}$. Then a $C^{2}$ function $f:\Upsilon \rightarrow {\Bbb %
R}$ is a convex function on the open convex set $\Upsilon $ if and only if
for each $x\in \Upsilon $, the Hessian of $f$ at $x$, denoted $({\rm Hess}%
f)_{x}$, is a positive semidefinite matrix.
\end{lem}

\begin{lem}
\label{lem-optimization-convex-2} {\rm (\cite[Corollary~21.2, p.~429]%
{Chong-Zak-01-book})} Let $\Upsilon \subset {\Bbb R}^{n}$ be an open convex
set in ${\Bbb R}^{n}$. Let $f:\Upsilon \rightarrow {\Bbb R}$ be a $C^{1}$
convex function with a point $x_{0}\in \Upsilon $ such that ${\rm grad}%
f\left( x_{0}\right) =0$, then the point $x_{0}$ is a global minimizer of $f$
over $\Upsilon $.
\end{lem}

For application purposes, we prove the following

\begin{lem}
\label{lem-convex-optimize-1} Let
\[
\Upsilon =\left\{ \left( x^{1},\ldots ,x^{n}\right) \in {\Bbb R}%
^{n}:x^{1}+\cdots +x^{n}=k\right\}
\]%
be a hyperplane of ${\Bbb R}^{n}$, and $f:{\Bbb R}^{n}\rightarrow {\Bbb R}$
a quadratic form given by
\begin{equation}
f\left( x^{1},\ldots ,x^{n}\right) =a\sum_{i=1}^{n-1}\left( x^{i}\right)
^{2}+b\left( x^{n}\right)^{2}-2\sum_{1\leq i<j\leq n}x^{i}x^{j},\qquad
a>0,\;b>0.  \label{eq-new-alg-f-1}
\end{equation}%
Then the constrained extremum problem
\begin{equation}
\min_{(x^{1},\ldots ,x^{n})\in \Upsilon}f  \label{eq-new-alg-f-2}
\end{equation}%
has a global solution given by
\begin{equation}
\left\{
\begin{array}{l}
\displaystyle x^{1}=x^{2}=\cdots =x^{n-1}=\frac{k}{a+1},\medskip \\
\displaystyle x^{n}=\frac{k}{b+1}=\frac{n-1}{b}\left( \frac{k}{a+1}\right)
=\left( a-n+2\right) \frac{k}{a+1},%
\end{array}
\right.  \label{eq-new-alg-f-3}
\end{equation}%
provided that
\begin{equation}
b=\frac{n-1}{a-n+2}.  \label{eq-rel-b-a-n-new}
\end{equation}
\end{lem}

\noindent {\bf Proof.} First we note that the set $\Upsilon $ is an open
convex set in ${\Bbb R}^{n}$ and the function $f$ is a $C^{\infty }$
function (and hence a $C^{2}$ function). Now we compute the matrix for the
Hessian ${\rm Hess}f$ of the function $f$. The partial derivatives of the
function $f$ are
\begin{equation}
\left\{
\begin{array}{l}
\displaystyle\frac{\partial f}{\partial x^{i}}=2\left( a+1\right)
x^{i}-2\sum_{\ell =1}^{n}x^{\ell },\quad i\in \{1,\ldots ,n-1\},\medskip \\
\displaystyle\frac{\partial f}{\partial x^{n}}=2\left( b+1\right)
x^{n}-2\sum_{\ell =1}^{n}x^{\ell }.%
\end{array}%
\right.  \label{eq-new-alg-f-4}
\end{equation}%
From (\ref{eq-new-alg-f-4}), we have
\begin{equation}
\left\{
\begin{array}{l}
\displaystyle\frac{\partial ^{2}f}{\partial \left( x^{i}\right) ^{2}}%
=2a,\quad i\in \{1,\ldots ,n-1\},\medskip \\
\displaystyle\frac{\partial ^{2}f}{\partial x^{i}\partial x^{j}}=-\,2,\quad
i,j\in \{1,\ldots ,n-1\},\medskip \\
\displaystyle\frac{\partial ^{2}f}{\partial x^{i}\partial x^{n}}=-\,2,\quad
i\in \{1,\ldots ,n-1\},\medskip \\
\displaystyle\frac{\partial ^{2}f}{\partial \left( x^{n}\right) ^{2}}=2b.%
\end{array}%
\right.  \label{eq-new-Hess-f-0}
\end{equation}%
Thus, in the standard frame of ${\Bbb R}^{n}$, the ${\rm Hess}f$ has the
matrix given by
\begin{equation}
2\left(
\begin{array}{ccccc}
a & -\,1 & \cdots & -\,1 & -\,1 \\
-\,1 & a & \cdots & -\,1 & -\,1 \\
\vdots & \vdots & \ddots & \vdots & \vdots \\
-\,1 & -\,1 & \cdots & a & -\,1 \\
-\,1 & -\,1 & \cdots & -\,1 & b%
\end{array}%
\right) .  \nonumber
\end{equation}%
We note that for any $X=\left( X^{1},\ldots ,X^{n}\right) \in T_{x}\Upsilon $%
, $x\in \Upsilon $, it follows that $\displaystyle\sum_{\ell =1}^{n}X^{\ell
}=0$. Consequently, for any $X=\left( X^{1},\ldots ,X^{n}\right) \in
T_{x}\Upsilon $, $x\in \Upsilon $ we have
\[
{\rm Hess}f\left( X,X\right) \geq 0.
\]%
Thus, for each $x\in \Upsilon $, the Hessian $({\rm Hess}f)_{x}$ of $f$ at $%
x $ is positive semidefinite. In view of Lemma~\ref%
{lem-optimization-convex-1}, this implies that the $C^{2}$ function $f$ is a
convex function on the open convex set $\Upsilon $.


For an optimal solution $\left( x^{1},\ldots ,x^{n}\right) $ of the problem (%
\ref{eq-new-alg-f-2}), the vector ${\rm grad}f$ is normal to $\Upsilon $,
equivalently, it is collinear with the vector $\left( 1,1,\ldots ,1\right) $%
. From (\ref{eq-new-alg-f-4}), for a critical point $x=\left( x^{1},\ldots
,x^{n}\right) $\ of the function $f$ we have
\begin{equation}
\left\{
\begin{array}{l}
\left( a+1\right) x^{i}-\sum_{\ell =1}^{n}x^{\ell }=0,\quad \quad i\in
\{1,\ldots ,n-1\},\medskip \\
\left( b+1\right) x^{n}-\sum_{\ell =1}^{n}x^{\ell }=0.%
\end{array}%
\right.  \label{eq-new-alg-f-4a}
\end{equation}%
From (\ref{eq-new-alg-f-4a}), it follows that a critical point $\left(
x^{1},\ldots ,x^{n-1},x^{n}\right) $ of the function $f$ has the form
\begin{equation}
x^{1}=\cdots =x^{n-1}=t,\quad x^{n}=\frac{n-1}{b}t.  \label{eq-new-alg-f-5}
\end{equation}%
Since
\[
x^{1}+x^{2}+\cdots +x^{n}=k,
\]%
in view of (\ref{eq-new-alg-f-5}), a critical point $\left( x^{1},\ldots
,x^{n}\right) $ of the considered problem is given by (\ref{eq-new-alg-f-3}%
). Solving one of the following three relations appearing in (\ref%
{eq-new-alg-f-3})
\[
\frac{k}{b+1}=\frac{n-1}{b}\left( \frac{k}{a+1}\right) =\left( a-n+2\right)
\frac{k}{a+1},
\]%
we get the equivalent relation given by (\ref{eq-rel-b-a-n-new}).
Consequently, in view of Lemma~\ref{lem-optimization-convex-2}, the point $%
\left( x^{1},\ldots ,x^{n}\right) $ given by (\ref{eq-new-alg-f-3}) is a
global minimum point. Inserting (\ref{eq-new-alg-f-3}) into (\ref%
{eq-new-alg-f-1}) we have $f\left( x^{1},\ldots ,x^{n}\right) =0$. $%
\blacksquare $

Now, we present the following Theorem, involving the Casorati inequalities
for algebraic Casorati curvatures $\delta _{{\cal C}^{T,\zeta }}\left(
r;n-1\right) $ and $\widehat{\delta }_{{\cal C}^{T,\zeta }}\left(
r;n-1\right) $.

\begin{th}
\label{th-alg-delta-C-(r;n-1)-1} Let $\left( M,g\right) $ be an {$n$%
-dimensional} Riemannian manifold, $\left( B,g_{B}\right) $ a Riemannian
vector bundle over $M$ and $\zeta $ a $B$-valued symmetric $\left(
1,2\right) $-tensor field. Let $T$ be a curvature-like tensor field
satisfying the algebraic Gauss equation {\rm (\ref{eq-Gauss-alg})}. Then
\begin{equation}
\left( \tau _{T}\right) _{{\rm Nor}}(p)\leq \frac{1}{n(n-1)}\left[ \delta _{%
{\cal C}^{T,\zeta }}(r;n-1) \right] _{p},\qquad 0<r<n\left( n-1\right) ,
\label{eq-alg-delta-C-(r;n-1)-ineq}
\end{equation}%
\begin{equation}
\left( \tau _{T}\right) _{{\rm Nor}}(p)\leq \frac{1}{n(n-1)}[\widehat{\delta
}_{{\cal C}^{T,\zeta }}(r;n-1) ]_{p},\qquad n(n-1) <r.
\label{eq-alg-delta-hat-C-(r;n-1)-ineq}
\end{equation}%
If
\begin{eqnarray*}
&\inf \{{\cal C}^{T,\zeta }(\Pi _{n-1}):\Pi _{n-1}\ {\rm is}\ {\rm a}\ {\rm %
hyperplane}\ {\rm of}\ T_{p}M\}& \\
&\left( {\rm resp.}\;\sup \{{\cal C}^{T,\zeta }(\Pi _{n-1}):\Pi _{n-1}\ {\rm %
is}\ {\rm a}\ {\rm hyperplane}\ {\rm of}\ T_{p}M\}\right) &
\end{eqnarray*}%
is attained by a hyperplane $\Pi _{n-1}$ of $T_{p}M$, $p\in M$, then the
equality sign holds in {\rm (\ref{eq-alg-delta-C-(r;n-1)-ineq})} {\rm (}%
resp. {\rm (\ref{eq-alg-delta-hat-C-(r;n-1)-ineq}))} if and only if with
respect to a suitable orthonormal tangent frame $\left\{
e_{1},...,e_{n}\right\} $ and a suitable orthonormal frame $\left\{
e_{n+1},...,e_{m}\right\} $ of the Riemann vector bundle $\left(
B,g_{B}\right) $, the components of $\zeta $ satisfy
\begin{equation}
\zeta _{ij}^{\alpha }=0\qquad i,j\in \{1,\ldots ,n\},\;\;i\neq j\;\;\alpha
\in \{n+1,\ldots ,m\},  \label{eq-alg-delta-C-(r;n-1)-zeta-1}
\end{equation}%
\begin{equation}
\zeta _{11}^{\alpha }=\zeta _{22}^{\alpha }=\cdots =\zeta
_{n-1\,n-1}^{\alpha }=\frac{r}{n(n-1) }\,\zeta _{nn}^{\alpha }\qquad \alpha
\in \{n+1,\ldots ,m\}.  \label{eq-alg-delta-C-(r;n-1)-zeta-2}
\end{equation}
\end{th}

\noindent {\bf Proof.} Let $p\in M$ and the set $\{e_{1},\ldots ,e_{n}\}$ be
an orthonormal basis of the tangent space $T_{p}M$ and $e_{\alpha }$ belong
to an orthonormal basis $\{e_{n+1},\ldots ,e_{m}\}$ of the Riemannian vector
bundle $(B,g_{B})$ over $M$ at $p$. We consider the following function
\begin{equation}
{\cal P}=r{\cal C}^{T,\zeta }+{\rm a}(r){\cal C}^{T,\zeta }(\Pi
_{n-1})-2\tau _{T}(p).  \label{eq-alg-Q-1}
\end{equation}%
where $\Pi _{n-1}$ is a hyperplane of $T_{p}M$. In view of (\ref%
{eq-alg-nC-trace-zeta}), the relation (\ref{eq-alg-Q-1}) becomes
\begin{equation}
{\cal P}=\left( n+r\right) {\cal C}^{T,\zeta }+{\rm a}(r){\cal C}^{T,\zeta
}(\Pi _{n-1})-\left\Vert {\rm trace}\,\zeta \right\Vert ^{2}.
\label{eq-alg-Q-2}
\end{equation}%
Without loss of generality, assume that the hyperplane $\Pi _{n-1}$ is
spanned by $e_{1},\ldots ,e_{n-1}$. Then from (\ref{eq-alg-Q-2}) it follows
that
\begin{equation}
{\cal P}=\frac{n+r}{n}\sum_{\alpha =n+1}^{m}\sum_{i,j=1}^{n}\left( \zeta
_{ij}^{\alpha }\right) ^{2}+\frac{{\rm a}(r)}{n-1}\sum_{\alpha
=n+1}^{m}\sum_{i,j=1}^{n-1}\left( \zeta _{ij}^{\alpha }\right)
^{2}-\sum_{\alpha =n+1}^{m}\left( \sum_{i=1}^{n}\zeta _{ii}^{\alpha }\right)
^{2}.  \label{eq-alg-Q-2a}
\end{equation}%
The function ${\cal P}$ is a quadratic polynomial in the components of the
tensor $\zeta $ and can be written as
\begin{eqnarray}
{\cal P} &=&\sum_{\alpha =n+1}^{m}\left\{ 2\left( \frac{r}{n}+\frac{{\rm a}%
(r)}{n-1}+1\right) \sum_{1\leq i<j\leq n-1}\left( \zeta _{ij}^{\alpha
}\right) ^{2}+\,2\left( \frac{r}{n}+1\right) \sum_{i=1}^{n-1}\left( \zeta
_{in}^{\alpha }\right) ^{2}\right.  \nonumber \\
&&\quad \left. +\,\left( \frac{r}{n}+\frac{{\rm a}(r)}{n-1}\right)
\sum_{i=1}^{n-1}\left( \zeta _{ii}^{\alpha }\right) ^{2}+\frac{r}{n}\left(
\zeta _{nn}^{\alpha }\right) ^{2}-2\sum_{1\leq i<j\leq n}\zeta _{ii}^{\alpha
}\zeta _{jj}^{\alpha }\right\}  \nonumber \\
&\geq &\sum_{\alpha =n+1}^{m}\left\{ \left( \frac{r}{n}+\frac{{\rm a}(r)}{n-1%
}\right) \sum_{i=1}^{n-1}\left( \zeta _{ii}^{\alpha }\right) ^{2}+\frac{r}{n}%
\left( \zeta _{nn}^{\alpha }\right) ^{2}-2\sum_{1\leq i<j\leq n}\zeta
_{ii}^{\alpha }\zeta _{jj}^{\alpha }\right\} .  \label{eq-alg-Q-3a}
\end{eqnarray}%
For $\alpha =n+1,\ldots ,m$, we consider a quadratic form
\[
f_{\alpha }:{\Bbb R}^{n}\rightarrow {\Bbb R}
\]%
given by
\begin{equation}
f_{\alpha }\left( \zeta _{11}^{\alpha },\ldots ,\zeta _{nn}^{\alpha }\right)
=\left( \frac{r}{n}+\frac{{\rm a}(r)}{n-1}\right) \sum_{i=1}^{n-1}\left(
\zeta _{ii}^{\alpha }\right) ^{2}+\frac{r}{n}\left( \zeta _{nn}^{\alpha
}\right) ^{2}-2\sum_{1\leq i<j\leq n}\zeta _{ii}^{\alpha }\zeta
_{jj}^{\alpha }  \label{eq-alg-f-r-alpha-0}
\end{equation}%
and the constrained extremum problem
\[
\min f_{\alpha },
\]%
subject to the condition
\[
\zeta _{11}^{\alpha }+\cdots +\zeta _{nn}^{\alpha }=k_{\alpha },
\]%
where $k_{\alpha }$ is a real constant. Comparing (\ref{eq-alg-f-r-alpha-0})
with (\ref{eq-new-alg-f-1}), we see that
\[
a=\left( \frac{r}{n}+\frac{{\rm a}(r)}{n-1}\right) ,\qquad b=\frac{r}{n},
\]%
which verifies the relation
\[
b=\frac{n-1}{a-n+2}
\]%
of (\ref{eq-rel-b-a-n-new}). Thus applying Lemma~\ref{lem-convex-optimize-1}%
, we see that the critical point
\[
\zeta ^{{\rm c}}=\left( \zeta _{11}^{\alpha },\zeta _{22}^{\alpha },\ldots
,\zeta _{n-1\,n-1}^{\alpha },\zeta _{nn}^{\alpha }\right)
\]%
given by
\begin{equation}
\zeta _{11}^{\alpha }=\zeta _{22}^{\alpha }=\cdots =\zeta
_{n-1\,n-1}^{\alpha }=\frac{r}{(n-1)\left( n+r\right) }\,k_{\alpha },\quad
\zeta _{nn}^{\alpha }=\frac{n}{n+r}\,k_{\alpha }  \label{eq-alg-f-r-alpha-1}
\end{equation}%
is a global minimum point. Inserting (\ref{eq-alg-f-r-alpha-1}) into (\ref%
{eq-alg-f-r-alpha-0}) we have $f_{\alpha }(\zeta ^{{\rm c}})=0$. Hence we
have
\begin{equation}
{\cal P}\geq 0,  \label{eq-Q-geq-0}
\end{equation}%
which in view of (\ref{eq-alg-Q-1}) gives
\begin{equation}
\frac{2\tau _{T}(p)}{n(n-1)}\leq \frac{r}{n(n-1)}\,{\cal C}_{p}^{T,\zeta }+%
\frac{{\rm a}(r)}{n(n-1)}\,{\cal C}^{T,\zeta }(\Pi _{n-1})
\label{eq-alg-Q-1a}
\end{equation}%
for every tangent hyperplane $\Pi _{n-1}$ of $T_{p}M$.

\smallskip

If $0<r<n(n-1)$, then ${\rm a}(r)>0$ and taking the infimum over all the
tangent hyperplanes $\Pi _{n-1}$ of $T_{p}M$, the relation (\ref{eq-alg-Q-1a}%
) gives the inequality (\ref{eq-alg-delta-C-(r;n-1)-ineq}). If $n(n-1)<r$,
then ${\rm a}(r)<0$, and taking the supremum over all the tangent
hyperplanes $\Pi _{n-1}$ of $T_{p}M$, the relation (\ref{eq-alg-Q-1a}) gives
the inequality (\ref{eq-alg-delta-hat-C-(r;n-1)-ineq}).

\smallskip

Suppose that
\begin{eqnarray*}
&\inf \{{\cal C}^{T,\zeta }(\Pi _{n-1}):\Pi _{n-1}\ {\rm is}\ {\rm a}\ {\rm %
hyperplane}\ {\rm of}\ T_{p}M\}& \\
&\left( {\rm resp.}\;\sup \{{\cal C}^{T,\zeta }(\Pi _{n-1}):\Pi _{n-1}\ {\rm %
is}\ {\rm a}\ {\rm hyperplane}\ {\rm of}\ T_{p}M\}\right) &
\end{eqnarray*}%
is attained by a hyperplane $\Pi _{n-1}$ spanned by $e_{1},\ldots ,e_{n-1}$.
Then the equality sign holds in (\ref{eq-alg-delta-C-(r;n-1)-ineq}) (resp. (%
\ref{eq-alg-delta-hat-C-(r;n-1)-ineq})) if and only if we have the equality
in all the previous inequalities. Thus the equality sign is true in the
inequality (\ref{eq-alg-delta-C-(r;n-1)-ineq}) (resp. (\ref%
{eq-alg-delta-hat-C-(r;n-1)-ineq})) if and only if the relations (\ref%
{eq-alg-delta-C-(r;n-1)-zeta-1}) and (\ref{eq-alg-delta-C-(r;n-1)-zeta-2})
are true. $\blacksquare $

Now, we have the following two results.

\begin{th}
\label{th-delta-C-(T,B)-(n-1)-1} Let $\left( M,g\right) $ be an {$n$%
-dimensional} Riemannian manifold, $\left( B,g_{B}\right) $ a Riemannian
vector bundle over $M$ and $\zeta $ a $B$-valued symmetric $\left(
1,2\right) $-tensor field. Let $T$ be a curvature-like tensor field
satisfying the algebraic Gauss equation {\rm (\ref{eq-Gauss-alg})}. Then the
$T$-normalized scalar curvature $(\tau _{T})_{{\rm Nor}}$ is bounded above
by the algebraic Casorati curvature $\delta _{{\cal C}^{T,\zeta }}(n-1)$
given by {\rm (\ref{eq-delta-C-(T,B)-(n-1)})}, that is,
\begin{equation}
(\tau _{T})_{{\rm Nor}}(p)\leq \left[ \delta _{{\cal C}^{T,\zeta }}(n-1)%
\right] _{p}.  \label{eq-delta-C-(T,B)-(n-1)-ineq}
\end{equation}%
If
\[
\inf \{{\cal C}^{T,\zeta }(\Pi _{n-1}):\Pi _{n-1}\ {\rm is}\ {\rm a}\ {\rm %
hyperplane}\ {\rm of}\ T_{p}M\}
\]%
is attained by a hyperplane $\Pi _{n-1}$ of $T_{p}M$, then the equality sign
holds in $(\ref{eq-delta-C-(T,B)-(n-1)-ineq})$ if and only if with respect
to suitable orthonormal tangent frame $\left\{ e_{1},...,e_{n}\right\} $ and
orthonormal frame $\left\{ e_{n+1},...,e_{m}\right\} $, the components of $%
\zeta $ satisfy
\begin{equation}
\zeta _{ij}^{\alpha }=0\qquad i,j\in \{1,\ldots ,n\},\;\;i\neq j\;\;\alpha
\in \{n+1,\ldots ,m\},  \label{eq-delta-C-(T,B)-(n-1)-zeta-1}
\end{equation}%
\begin{equation}
\zeta _{11}^{\alpha }=\zeta _{22}^{\alpha }=\cdots =\zeta
_{n-1\,n-1}^{\alpha }=\frac{1}{2}\,\zeta _{nn}^{\alpha }\,,\qquad \alpha \in
\{n+1,\ldots ,m\}.  \label{eq-delta-C-(T,B)-(n-1)-zeta-2}
\end{equation}
\end{th}

\noindent {\bf Proof.} Using
\begin{equation}
\left[ \delta _{{\cal C}^{T,\zeta }}(n-1)\right] _{p}=\frac{1}{n(n-1)}\left[
\delta _{{\cal C}^{T,\zeta }}\left( \frac{n(n-1)}{2};n-1\right) \right] _{p}
\label{eq-delta-C-(T,B)-(r;n-1)-C-(T,B)-(n-1)}
\end{equation}%
in (\ref{eq-alg-delta-C-(r;n-1)-ineq}), we get (\ref%
{eq-delta-C-(T,B)-(n-1)-ineq}). Taking $2r=n(n-1)$ in (\ref%
{eq-alg-delta-C-(r;n-1)-zeta-2}) we get (\ref{eq-delta-C-(T,B)-(n-1)-zeta-2}%
). $\blacksquare $

\begin{th}
\label{th-delta-hat-C-(T,B)-(n-1)-1} Let $\left( M,g\right) $ be an {$n$%
-dimensional} Riemannian manifold, $\left( B,g_{B}\right) $ a Riemannian
vector bundle over $M$ and $\zeta $ a $B$-valued symmetric $\left(
1,2\right) $-tensor field. Let $T$ be a curvature-like tensor field
satisfying the algebraic Gauss equation {\rm (\ref{eq-Gauss-alg})}. Then the
$T$-normalized scalar curvature $(\tau _{T})_{{\rm Nor}}$ is bounded above
by the algebraic Casorati curvature $\widehat{\delta }_{{\cal C}^{T,\zeta
}}(n-1)$, that is,
\begin{equation}
(\tau _{T})_{{\rm Nor}}(p)\leq \lbrack \widehat{\delta }_{{\cal C}^{T,\zeta
}}(n-1)]_{p}.  \label{eq-delta-hat-C-(T,B)-(n-1)-ineq}
\end{equation}%
If
\[
\sup \{{\cal C}^{T,\zeta }(\Pi _{n-1}):\Pi _{n-1}\ {\rm is}\ {\rm a}\ {\rm %
hyperplane}\ {\rm of}\ T_{p}M\}
\]%
is attained by a hyperplane $\Pi _{n-1}$ of $T_{p}M$, then the equality sign
in {\rm (\ref{eq-delta-hat-C-(T,B)-(n-1)-ineq})} is true if and only if with
respect to a suitable orthonormal tangent frame $\left\{
e_{1},...,e_{n}\right\} $ and a suitable orthonormal frame $\left\{
e_{n+1},...,e_{m}\right\} $ of the Riemann vector bundle $\left(
B,g_{B}\right) $, the components of $\zeta $ satisfy
\begin{equation}
\zeta _{ij}^{\alpha }=0,\qquad i,j\in \{1,\ldots ,n\},\quad i\neq j,\quad
\alpha \in \{n+1,\ldots ,m\},  \label{eq-delta-hat-C-(T,B)-(n-1)-zeta-1}
\end{equation}%
\begin{equation}
\zeta _{11}^{\alpha }=\zeta _{22}^{\alpha }=\cdots =\zeta
_{n-1\,n-1}^{\alpha }=2\,\zeta _{nn}^{\alpha }\,,\qquad \alpha \in
\{n+1,\ldots ,m\}.  \label{eq-delta-hat-C-(T,B)-(n-1)-zeta-2}
\end{equation}
\end{th}

\noindent {\bf Proof.} Using
\begin{equation}
\left[ \widehat{\delta }_{{\cal C}^{T,\zeta }}(n-1)\right] _{p}=\frac{1}{%
n(n-1)}\left[ \widehat{\delta }_{{\cal C}^{T,\zeta }}\left(
2n(n-1);n-1\right) \right] _{p}
\label{eq-delta-hat-C-(T,B)-(r;n-1)-C-(T,B)-(n-1)}
\end{equation}%
in (\ref{eq-alg-delta-hat-C-(r;n-1)-ineq}), we get (\ref%
{eq-delta-hat-C-(T,B)-(n-1)-ineq}). Taking $r=2n(n-1)$ in (\ref%
{eq-alg-delta-C-(r;n-1)-zeta-2}) we get (\ref%
{eq-delta-hat-C-(T,B)-(n-1)-zeta-2}). $\blacksquare $

\section{Casorati inequalities for Riemannian submanifolds\label%
{sect-Casorati-ineq-Riem-submfd}}

\begin{th}
\label{th-delta-C-(r;n-1)-1} Let $(M,g)$ be an {$n$-dimensional} Riemannian
submanifold of {$m$-dimensional} Riemannian manifold $(\widetilde{M},%
\widetilde{g})$. Then the generalized normalized {$\delta $-Casorati}
curvatures $\delta _{{\cal C}}(r;n-1) $ and $\widehat{\delta }_{{\cal C}%
}(r;n-1) $ satisfy
\begin{equation}
\tau _{{\rm Nor}}(p)\leq \frac{1}{n(n-1)}\,[\delta _{{\cal C}}\left(
r;n-1\right) ]_{p}+\widetilde{\tau }_{{\rm Nor}}\left( T_{p}M\right) ,\quad
0<r<n(n-1) ,  \label{eq-delta-C-(r;n-1)-ineq}
\end{equation}%
\begin{equation}
\tau _{{\rm Nor}}(p)\leq \frac{1}{n(n-1)}\,[\widehat{\delta }_{{\cal C}%
}(r;n-1) ]_{p}+\widetilde{\tau }_{{\rm Nor}}\left( T_{p}M\right) ,\quad
n(n-1) <r.  \label{eq-delta-hat-C-(r;n-1)-ineq}
\end{equation}%
If \vspace{-5mm}
\begin{eqnarray*}
&\inf \{{\cal C}(\Pi _{n-1}):\Pi _{n-1}\ {\rm is}\ {\rm a}\ {\rm hyperplane}%
\ {\rm of}\ T_{p}M\}& \\
&\left( {\rm resp.}\;\sup \{{\cal C}(\Pi _{n-1}):\Pi _{n-1}\ {\rm is}\ {\rm a%
}\ {\rm hyperplane}\ {\rm of}\ T_{p}M\}\right) &
\end{eqnarray*}%
is attained by a hyperplane $\Pi _{n-1}$ of $T_{p}M$, $p\in M$, then the
equality sign holds in {\rm (\ref{eq-delta-C-(r;n-1)-ineq})} {\rm (}resp.
{\rm (\ref{eq-delta-hat-C-(r;n-1)-ineq}))} for all $p\in M$ if and only if $%
\left( M,g\right) $ is an invariantly quasi-umbilical submanifold with
trivial normal connection in $(\widetilde{M},\widetilde{g})$, such that with
respect to suitable tangent orthonormal frame $\left\{
e_{1},...,e_{n}\right\} $ and normal orthonormal frame $\left\{
e_{n+1},...,e_{m}\right\} $, the shape operators $A_{\alpha }\equiv
A_{e_{\alpha }}$, $\alpha \in \{n+1,...,m\}$, take the following forms:
\begin{equation}
A_{n+1}=\left(
\begin{array}{cccccc}
a & 0 & 0 & ... & 0 & 0 \\
0 & a & 0 & ... & 0 & 0 \\
0 & 0 & a & ... & 0 & 0 \\
\vdots & \vdots & \vdots & \ddots & \vdots & \vdots \\
0 & 0 & 0 & ... & a & 0 \\
0 & 0 & 0 & ... & 0 & \displaystyle\frac{n(n-1) }{r}\,a%
\end{array}%
\right) ,\ A_{n+2}=\cdots =A_{m}=0.  \label{eq-delta-C-(r;n-1)-shape}
\end{equation}
\end{th}

\noindent {\bf Proof.} Let $\left( M,g\right) $ be an {$n$-dimensional}
Riemannian submanifold of an {$m$-dimensional} Riemannian manifold $(%
\widetilde{M},\widetilde{g})$. Let the Riemannian vector bundle $\left(
B,g_{B}\right) $ over $M$ be replaced by the normal bundle $T^{\perp }M$,
and the $B$-valued symmetric $\left( 1,2\right) $-tensor field $\zeta $ be
replaced by the second fundamental form of immersion $\sigma $. In (\ref%
{eq-Gauss-alg}), we set
\[
T\left( X,Y,Z,W\right) =R(X,Y,Z,W)-\widetilde{R}(X,Y,Z,W)
\]%
with $R$ the Riemann curvature tensor on $M$ and $\zeta =\sigma $. Then we
see that
\[
(\tau _{T})_{{\rm Nor}}(p)=\tau _{{\rm Nor}}(p)-\widetilde{\tau }_{{\rm Nor}%
}\left( T_{p}M\right),
\]%
\[
\delta _{{\cal C}^{T,\zeta }}(r;n-1) =\delta _{{\cal C}}\left( r;n-1\right)
,
\]
\[
\widehat{\delta }_{{\cal C}^{T,\zeta }}(r;n-1) =\widehat{\delta }_{{\cal C}%
}(r;n-1) .
\]%
Using these facts in (\ref{eq-alg-delta-C-(r;n-1)-ineq}) and (\ref%
{eq-alg-delta-hat-C-(r;n-1)-ineq}), we get (\ref{eq-delta-C-(r;n-1)-ineq})
and (\ref{eq-delta-hat-C-(r;n-1)-ineq}), respectively.

\smallskip

The conditions of equality cases (\ref{eq-alg-delta-C-(r;n-1)-zeta-1}) and (%
\ref{eq-alg-delta-C-(r;n-1)-zeta-2}) become
\begin{equation}
\sigma _{ij}^{\alpha }=0\qquad i,j\in \{1,\ldots ,n\},\;\;i\neq j\;\;\alpha
\in \{n+1,\ldots ,m\}  \label{eq-delta-C-(r;n-1)-shape-1}
\end{equation}%
and
\begin{equation}
\sigma _{11}^{\alpha }=\sigma _{22}^{\alpha }=\cdots =\sigma
_{n-1\,n-1}^{\alpha }=\frac{r}{n(n-1) }\,\sigma _{nn}^{\alpha },\qquad
\alpha \in \{n+1,\ldots ,m\},  \label{eq-delta-C-(r;n-1)-shape-2}
\end{equation}%
respectively. Thus the equality sign holds in both the inequalities (\ref%
{eq-delta-C-(r;n-1)-ineq}) and (\ref{eq-delta-hat-C-(r;n-1)-ineq}) if and
only if (\ref{eq-delta-C-(r;n-1)-shape-1}) and (\ref%
{eq-delta-C-(r;n-1)-shape-2}) are true.

The interpretation of the relations (\ref{eq-delta-C-(r;n-1)-shape-1}) is
that the shape operators with respect to all normal directions $e_{\alpha }$
commute, or equivalently, that the normal connection $\nabla ^{\perp }$ is
flat, or still, that the {\em normal curvature tensor} $R^{\perp }$, that
is, the curvature tensor of the normal connection, is trivial. Furthermore,
the interpretation of the relations (\ref{eq-delta-C-(r;n-1)-shape-2}) is
that there exist $m-n$ mutually orthogonal unit normal vectors $\left\{
e_{n+1},...,e_{m}\right\} $ such that the shape operators with respect to
all directions $e_{\alpha }$ $\left( \alpha \in \left\{
e_{n+1},...,e_{m}\right\} \right) $ have an eigenvalue of multiplicity $n-1$
and that for each $e_{\alpha }$ the distinguished eigendirection is the same
(namely $e_{n}$), that is, the submanifold is {\em invariantly
quasi-umbilical} \cite{Bl-Led-77-Stevin}.

Thus from the relations (\ref{eq-delta-C-(r;n-1)-shape-1}) and (\ref%
{eq-delta-C-(r;n-1)-shape-2}), we conclude that the equality holds in (\ref%
{eq-delta-C-(r;n-1)-ineq}) and/or (\ref{eq-delta-hat-C-(r;n-1)-ineq}) for
all $p\in M$ if and only if the Riemannian submanifold $M$ is invariantly
quasi-umbilical with trivial normal connection $\nabla ^{\perp }$ in $%
\widetilde{M}$, such that with respect to suitable orthonormal tangent and
normal orthonormal frames, the shape operators take the form given by (\ref%
{eq-delta-C-(r;n-1)-shape}). $\blacksquare $

\begin{th}
\label{th-delta-C-(n-1)-1} Let $\left( M,g\right) $ be an {$n$-dimensional}
Riemannian submanifold of {$m$-dimensional} Riemannian manifold $(\widetilde{%
M},\widetilde{g})$. Then the normalized {$\delta $-Casorati} curvature $%
\delta _{{\cal C}}(n-1)$ satisfies
\begin{equation}
\tau _{{\rm Nor}}(p)\leq \lbrack \delta _{{\cal C}}(n-1)]_{p}+\widetilde{%
\tau }_{{\rm Nor}}\left( T_{p}M\right) .  \label{eq-delta-C-(n-1)-ineq}
\end{equation}%
If \vspace{-5mm}
\[
\inf \{{\cal C}(\Pi _{n-1}):\Pi _{n-1}\ {\rm is}\ {\rm a}\ {\rm hyperplane}\
{\rm of}\ T_{p}M\}
\]%
is attained by a hyperplane $\Pi _{n-1}$ of $T_{p}M$, $p\in M$, then the
equality sign holds if and only if $M$ is an invariantly quasi-umbilical
submanifold with trivial normal connection in $\widetilde{M}$, such that
with respect to suitable orthonormal tangent frame $\left\{
e_{1},...,e_{n}\right\} $ and normal orthonormal frame $\left\{
e_{n+1},...,e_{m}\right\} $, the shape operators $A_{\alpha }\equiv
A_{e_{\alpha }}$, $\alpha \in \{n+1,...,m\}$, take the following forms
\begin{equation}
A_{n+1}=\left(
\begin{array}{cccccc}
a & 0 & 0 & ... & 0 & 0 \\
0 & a & 0 & ... & 0 & 0 \\
0 & 0 & a & ... & 0 & 0 \\
\vdots & \vdots & \vdots & \ddots & \vdots & \vdots \\
0 & 0 & 0 & ... & a & 0 \\
0 & 0 & 0 & ... & 0 & 2a%
\end{array}%
\right) ,\ \quad A_{n+2}=\cdots =A_{m}=0.  \label{eq-delta-C-(n-1)-shape}
\end{equation}
\end{th}

\noindent {\bf Proof.} Using (\ref{eq-delta-C-(r;n-1)-C-(n-1)}) in (\ref%
{eq-delta-C-(r;n-1)-ineq}), we get (\ref{eq-delta-C-(n-1)-ineq}). Putting $%
2r=n(n-1)$ in (\ref{eq-delta-C-(r;n-1)-shape}) we get (\ref%
{eq-delta-C-(n-1)-shape}). $\blacksquare $

\begin{th}
\label{th-delta-hat-C-(n-1)-1} Let $\left( M,g\right) $ be an {$n$%
-dimensional} Riemannian submanifold of {$m$-dimensional} Riemannian
manifold $(\widetilde{M},\widetilde{g})$. Then the normalized {$\delta $%
-Casorati} curvature $\widehat{\delta }_{{\cal C}}(n-1)$ satisfies
\begin{equation}
\tau _{{\rm Nor}}(p)\leq \lbrack \widehat{\delta }_{{\cal C}}(n-1)]_{p}+%
\widetilde{\tau }_{{\rm Nor}}\left( T_{p}M\right) .
\label{eq-delta-hat-C-(n-1)-ineq}
\end{equation}%
If \vspace{-5mm}
\[
\sup \{{\cal C}(\Pi _{n-1}):\Pi _{n-1}\ {\rm is}\ {\rm a}\ {\rm hyperplane}\
{\rm of}\ T_{p}M\}
\]%
is attained by a hyperplane $\Pi _{n-1}$ of $T_{p}M$, $p\in M$, then the
equality sign holds if and only if $\left( M,g\right) $ is an invariantly
quasi-umbilical submanifold with trivial normal connection in $(\widetilde{M}%
,\widetilde{g})$, such that with respect to suitable orthonormal tangent
frame $\left\{ e_{1},...,e_{n}\right\} $ and normal orthonormal frame $%
\left\{ e_{n+1},...,e_{m}\right\} $, the shape operators $A_{\alpha }\equiv
A_{e_{\alpha }}$, $\alpha \in \{n+1,...,m\}$, take the following forms:
\begin{equation}
A_{n+1}=\left(
\begin{array}{cccccc}
a & 0 & 0 & ... & 0 & 0 \\
0 & a & 0 & ... & 0 & 0 \\
0 & 0 & a & ... & 0 & 0 \\
\vdots & \vdots & \vdots & \ddots & \vdots & \vdots \\
0 & 0 & 0 & ... & a & 0 \\
0 & 0 & 0 & ... & 0 & \displaystyle\frac{1}{2}\,a%
\end{array}%
\right) ,\ A_{n+2}=\cdots =A_{m}=0.  \label{eq-delta-hat-C-(n-1)-shape}
\end{equation}
\end{th}

\noindent {\bf Proof.} Using (\ref{eq-delta-hat-C-(r;n-1)-C-(n-1)}) in (\ref%
{eq-delta-hat-C-(r;n-1)-ineq}), we get (\ref{eq-delta-hat-C-(n-1)-ineq}).
Putting $r=2n(n-1)$ in (\ref{eq-delta-C-(r;n-1)-shape}) we get (\ref%
{eq-delta-hat-C-(n-1)-shape}). $\blacksquare $

\section{Casorati inequalities for submanifolds of real space forms\label%
{sect-Casorati-ineq-Riem-submfd-RSF}}

An {$m$-dimensional} Riemannian manifold $(\widetilde{M},\widetilde{g})$
with constant sectional curvature $c$, denoted $\widetilde{M}(c)$, is called
a {\em real space form}, and its Riemann curvature tensor $\widetilde{R}$ is
then given by
\begin{equation}
\widetilde{R}(X,Y,Z,W)=c\left\{ \widetilde{g}\left( Y,Z\right) \widetilde{g}%
\left( X,W\right) -\widetilde{g}\left( X,Z\right) \widetilde{g}\left(
Y,W\right) \right\}  \label{eq-curv-RSF}
\end{equation}%
for all vector fields $X,Y,Z,W$ on $\widetilde{M}$. The model spaces for
real space forms are the Euclidean spaces ($c=0$), the spheres ($c>0$), and
the hyperbolic spaces ($c<0$). For an {$n$-dimensional} Riemannian
submanifold $\left( M,g\right) $ of a real space form $\widetilde{M}(c)$ it
is easy to see that
\begin{equation}
\widetilde{\tau }_{{\rm Nor}}\left( T_{p}M\right) =c.
\label{eq-tau-tilde-Nor-RSF}
\end{equation}

\begin{th}
\label{th-delta-C-(r;n-1)-RSF} {\rm \cite[Theorem~2.1 and Corollary~3.1]%
{Decu-HV-08-JIPAM}} Let $\left( M,g\right) $ be an {$n$-dimensional}
Riemannian submanifold of {$m$-dimensional} real space form $\widetilde{M}%
(c) $. Then
\begin{equation}
\tau _{{\rm Nor}}(p)\leq \frac{1}{n(n-1) }\,[\delta _{{\cal C}}(r;n-1)
]_{p}+c,\quad 0<r<n(n-1) ,  \label{eq-delta-C-(r;n-1)-ineq-RSF}
\end{equation}%
\begin{equation}
\tau _{{\rm Nor}}(p)\leq \frac{1}{n(n-1) }\,[\widehat{\delta }_{{\cal C}%
}(r;n-1) ]_{p}+c,\quad n(n-1) <r.  \label{eq-delta-hat-C-(r;n-1)-ineq-RSF}
\end{equation}%
The equality sign holds in {\rm (\ref{eq-delta-C-(r;n-1)-ineq-RSF})} {\rm (}%
resp. {\rm (\ref{eq-delta-hat-C-(r;n-1)-ineq-RSF}))} for all $p\in M$ if and
only if $\left( M,g\right) $ is an invariantly quasi-umbilical submanifold
with trivial normal connection in $\widetilde{M}(c)$, such that with respect
to suitable tangent orthonormal frame $\left\{ e_{1},...,e_{n}\right\} $ and
normal orthonormal frame $\left\{ e_{n+1},...,e_{m}\right\} $, the shape
operators $A_{\alpha }\equiv A_{e_{\alpha }}$, $\alpha \in \{n+1,...,m\}$,
take the forms given by $(\ref{eq-delta-C-(r;n-1)-shape})$.
\end{th}

\noindent {\bf Proof.} Using (\ref{eq-tau-tilde-Nor-RSF}) in (\ref%
{eq-delta-C-(r;n-1)-ineq}) and (\ref{eq-delta-hat-C-(r;n-1)-ineq}) we get (%
\ref{eq-delta-C-(r;n-1)-ineq-RSF}) and (\ref{eq-delta-hat-C-(r;n-1)-ineq-RSF}%
), respectively. $\blacksquare $

\begin{th}
\label{th-delta-C-(n-1)-RSF} {\rm (Theorem~4.1, \cite{ZhangP-Z-15-arXiv})}
Let $\left( M,g\right) $ be an {$n$-dimensional} Riemannian submanifold of {$%
m$-dimensional} real space form $\widetilde{M}(c) $. Then the normalized {$%
\delta $-Casorati} curvature $\delta _{{\cal C}}(n-1)$ satisfies
\begin{equation}
\tau _{{\rm Nor}}(p)\leq \lbrack \delta _{{\cal C}}(n-1)]_{p}+c.
\label{eq-delta-C-(n-1)-ineq-RSF}
\end{equation}%
Moreover, the equality sign holds for all $p\in M$ if and only if $\left(
M,g\right) $ is an invariantly quasi-umbilical submanifold with trivial
normal connection in $(\widetilde{M},\widetilde{g})$, such that with respect
to suitable orthonormal tangent frame $\left\{ e_{1},...,e_{n}\right\} $ and
normal orthonormal frame $\left\{ e_{n+1},...,e_{m}\right\} $, the shape
operators $A_{\alpha }\equiv A_{e_{\alpha }}$, $\alpha \in \{n+1,...,m\}$,
take the forms given by $(\ref{eq-delta-C-(n-1)-shape})$.
\end{th}

\noindent {\bf Proof.} Using (\ref{eq-delta-C-(r;n-1)-C-(n-1)}) in (\ref%
{eq-delta-C-(r;n-1)-ineq-RSF}), we get (\ref{eq-delta-C-(n-1)-ineq-RSF}). $%
\blacksquare $

\begin{th}
\label{th-delta-hat-C-(n-1)-RSF} {\rm (Theorem~1 and Corollary~3, \cite%
{Decu-HV-07-Brasov})} Let $\left( M,g\right) $ be an {$n$-dimensional}
Riemannian submanifold of {$m$-dimensional} real space form $\widetilde{M}%
(c) $. Then the normalized {$\delta $-Casorati} curvature $\widehat{\delta }%
_{{\cal C}}(n-1)$ satisfies
\begin{equation}
\tau _{{\rm Nor}}(p)\leq \lbrack \widehat{\delta }_{{\cal C}}(n-1)]_{p}+c.
\label{eq-delta-hat-C-(n-1)-ineq-RSF}
\end{equation}%
Moreover, the equality sign holds for all $p\in M$ if and only if $\left(
M,g\right) $ is an invariantly quasi-umbilical submanifold with trivial
normal connection in $(\widetilde{M},\widetilde{g})$, such that with respect
to suitable orthonormal tangent frame $\left\{ e_{1},...,e_{n}\right\} $ and
normal orthonormal frame $\left\{ e_{n+1},...,e_{m}\right\} $, the shape
operators $A_{\alpha }\equiv A_{e_{\alpha }}$, $\alpha \in \{n+1,...,m\}$,
take the forms given by $(\ref{eq-delta-hat-C-(n-1)-shape})$.
\end{th}

\noindent {\bf Proof.} Using (\ref{eq-delta-hat-C-(r;n-1)-C-(n-1)}) in (\ref%
{eq-delta-hat-C-(r;n-1)-ineq-RSF}), we get (\ref%
{eq-delta-hat-C-(n-1)-ineq-RSF}). $\blacksquare $

\section{Further studies\label{sect-Casorati-further-studies}}

In this section, we present some problems. Similar problems can be
formulated in those situations, where Riemman curvature tensor of the
ambient manifold has some nice well known form.

\begin{prob-new}
Like in \cite{Chen-05-AMH}, to obtain Casorati inequalities for conformally
flat submanifolds of a real space form.
\end{prob-new}

\begin{prob-new}
Riemannian manifolds of quasi-constant curvature (cf.\ \cite{Boju-Pop-78-JDG}%
, \cite{Chen-Yano-72}, \cite{Gan-Mih-00-Angew}, \cite{Mocanu-86}, \cite%
{Vranceanu-68}) represent a good generalization of real space forms. To
obtain Casorati inequalities for submanifolds of quasi-constant curvature
manifolds. To study Casorati ideal submanifolds of quasi-constant curvature
manifolds.
\end{prob-new}

\begin{prob-new}
To obtain Casorati inequalities for submanifolds of generalized complex
space forms (cf.\ \cite{KimUK-99}, \cite{Olszak-89}, \cite{Vanhecke-75/76},
\cite{TV-81-TAMS}).
\end{prob-new}

\begin{prob-new}
Like the improved Chen-Ricci inequalities \cite{Tri-11-DGA}, to improve
Casorati inequalities for Lagrangian \cite{Chen-01-Taiwan} and Kaehlerian
slant submanifolds \cite{Chen-90-slant} of a complex space form, if possible.
\end{prob-new}

\begin{prob-new}
To obtain Casorati inequalities for different kind of submanifolds of
locally conformal Kaehler space forms (cf.\ \cite{Drag-Ornea-99-book}, \cite%
{Vaisman-76}).
\end{prob-new}

\begin{prob-new}
Like the improved Chen-Ricci inequalities \cite{Tri-11-DGA}, to improve
Casorati inequalities for Lagrangian submanifolds of a locally conformal
Kaehler space form (under some conditions), if possible.
\end{prob-new}

\begin{prob-new}
To obtain Casorati inequalities for submanifolds of Kaehler manifolds of
quasi constant holomorphic sectional curvatures (cf. \cite{Gan-Mih-08-CEJM},
\cite{Bej-Ben-08-JG}).
\end{prob-new}

\begin{prob-new}
To obtain Casorati inequalities for different kind of submanifolds of
Bochner-Kaehler manifolds \cite{Chen-Yano-75}.
\end{prob-new}

\begin{prob-new}
Like the improved Chen-Ricci inequalities \cite{Tri-11-DGA}, to improve
Casorati inequalities for Lagrangian submanifolds of Bochner-Kaehler
manifolds, if possible.
\end{prob-new}

\begin{prob-new}
Like the improved Chen-Ricci inequalities \cite{Tri-11-DGA}, to improve
Casorati inequalities for Lagrangian submanifolds of a quaternionic space
form \cite{Ishihara-74}, if possible.
\end{prob-new}

\begin{prob-new}
To obtain Casorati inequalities for different kind of submanifolds \cite%
{Tri-96-JIMS} of generalized $\left( \kappa ,\mu \right) $ space forms \cite%
{Car-Mart-Tri-13-MedJM} and in particular generalized Sasakian space forms
\cite{Aleg-BC-04-Israel} and Sasakian space forms.
\end{prob-new}

\begin{prob-new}
Like the improved Chen-Ricci inequalities \cite{Tri-11-DGA}, to improve
Casorati inequalities for Legendrian submanifolds of a Sasakian space form
(cf. \cite{Tanno-69-Tohoku-2}, \cite{Bl-10-book}).
\end{prob-new}

\begin{prob-new}
To obtain Casorati inequalities for different kind of submanifolds of
different kind of manifolds equipped with a semi-symmetric metric connection
(cf. \cite{Pak-69}, \cite{Yano-70-semi}, \cite{Nakao-76}).
\end{prob-new}

\begin{prob-new}
To obtain Casorati inequalities for centroaffine hypersurfaces \cite%
{Nom-Sasaki-94-Aff-Diff-Geom}.
\end{prob-new}

\noindent {\bf Acknowledgements.} The author is thankful to
Professor Ugo Gianazza (gianazza@imati.cnr.it), Claudio Gnoli
(claudio.gnoli@unipv.it), Claudia Olivati and Anna Bendiscioli from
University of Pavia, Italy for their help in tracing the original
papers of Felice Casorati (\cite{Casorati-1889},
\cite{Casorati-1889-note}).

\noindent Department of Mathematics

\noindent Institute of Science

\noindent Banaras Hindu University

\noindent Varanasi 221005, India

\noindent Email: mmtripathi66@yahoo.com

\end{document}